\documentclass[9pt,journal,twoside,web]{ieeecolor}
\usepackage{generic}
\usepackage{cite}
\usepackage{amsmath,amssymb,amsfonts}
\usepackage{algorithmic}
\usepackage{graphicx}
\usepackage{algorithm,algorithmic}
\usepackage{hyperref}
\hypersetup{hidelinks=true}
\usepackage{textcomp}
\usepackage{CJKutf8}

\newtheorem{Remark}{\bf Remark}[section]

\newtheorem{Definition}{\bf Definition}[section]

\newtheorem{Problem}{\bf Problem}[section]

\newenvironment{Proof}{\noindent{\em Proof:\/}}{\hfill $\Box$\par}

\newtheorem{Theorem}{\bf Theorem}[section]

\newtheorem{Lemma}{\bf Lemma}[section]

\def\BibTeX{{\rm B\kern-.05em{\sc i\kern-.025em b}\kern-.08em
    T\kern-.1667em\lower.7ex\hbox{E}\kern-.125emX}}
\begin{document}
\begin{CJK}{UTF8}{gbsn}
\title{Attitude Tracking of Uncertain Flexible Spacecraft Systems Subject to Unknown External Disturbances}
\author{Zean Bao, Maobin Lu, \IEEEmembership{Member, IEEE}, Fang Deng,  \IEEEmembership{Senior Member, IEEE}, and Jie Chen, \IEEEmembership{Fellow, IEEE}
\thanks{This work was supported
in part by National Key R\&D Program of China under Grant 2021ZD0112600,
in part by National Natural Science Foundation of China under Grant
62373058,
in part by Key Program of National Natural Science Foundation of China under Grant 61933002 and U20B2073,
in part by National Science Fund for Distinguished Young Scholars of China under Grant 62025301,
in part by Natural Science Foundation of Chongqing under Grant 2021ZX4100036, and
in part by the Basic Science Center Programs of NSFC under Grant 62088101.}
\thanks{Zean Bao, Maobin Lu, and Fang Deng are with the School of
Automation, Beijing Institute of Technology, Beijing 100081, China, and
also with the Beijing Institute of Technology Chongqing Innovation
Center, Chongqing 401135, China (e-mail: baozean@gmail.com;
lumaobin@bit.edu.cn; dengfang@bit.edu.cn) }
\thanks{Jie Chen is with the Beijing Institute of Technology, Beijing 100081, China, and also with Tongji University, Shanghai 200092, China (E-mail: chenjie@bit.edu.cn)}
}

\maketitle

\begin{abstract}
In this paper, we investigate the attitude tracking problem of uncertain flexible spacecraft systems subject to external disturbances.
In sharp contrast to existing results, the dynamics of flexible spacecraft systems and external disturbances are allowed to be unknown.
To deal with the challenges by these unknown factors, we develop a class of nonlinear internal models which converts
the attitude tracking problem of uncertain flexible spacecraft systems 
into a regulation problem of an augmented system.
Furthermore, to overcome the difficulties caused by the unmeasurable modal variable, the uncertainty introduced by the internal model, and the cross-coupling of the uncertainties with the system state, 
we design an auxiliary dynamic system for auxiliary stabilization, a dynamic compensator for dynamic compensation, and a linearly parameterized transformation for adaptive regulation in sequence. By introducing a series of coordinate and input transformations, we propose an adaptive dynamic control law to achieve regulation of the augmented system and thus leading to the solution to the attitude tracking problem.
In addition, we analyze the convergence issue of the estimated parameter to its true value by the persistently exciting condition. 
Finally, the effectiveness of the developed approach is verified by its application to the attitude manoeuvre of a flexible spacecraft system in the presence of external disturbances.
\end{abstract}

\begin{IEEEkeywords}
Uncertain flexible spacecraft system, tracking control, adaptive control, external disturbance. 
\end{IEEEkeywords}

\section{Introduction}

\IEEEPARstart{A}{ttitude} control is essential for maintaining stability of the spacecraft. 
With increasingly complex work tasks, accurate attitude control of the spacecraft is becoming more and more important and thus, the attitude control problem has been considered by a lot of researchers (cf., \cite{kelkar1995dissipative,di2002output,chen2009attitude,chen2014attitude,zhong2017attitude,ding2014nonsmooth,cao2019faster,yu2021antidisturbance,liu2021adaptive,zhang2021fixed,ahmed1998adaptive,kang1995nonlinear,dong2022anti,meng2023second}). 



Over the past few decades, significant endeavours have been dedicated to attitude control of rigid spacecraft systems (cf., \cite{joshi1995robust,ahmed1998adaptive,mayhew2011quaternion,benallegue2018adaptive,tayebi2008unit,hegrenaes2005spacecraft,lu2018attitude,lu2019leader,zou2020fixed}). In particular, a nonlinear adaptive control law is designed for the attitude control problem of uncertain rigid spacecraft systems \cite{ahmed1998adaptive} and the global asymptotic attitude tracking is obtained.  
A unit quaternion-based dynamic output feedback control law without angular velocity measurement is developed in \cite{tayebi2008unit} to solve the attitude tracking problem of rigid spacecraft systems. 
A quaternion-based hybrid feedback scheme is further proposed to study the similar problem of the rigid spacecraft systems in \cite{mayhew2011quaternion}. This approach can avoid unwinding
and eliminate extreme measurement noise sensitivity in \cite{mayhew2011quaternion}. A sliding mode control law based on the modified Rodrigues parameters is developed in \cite{dong2021anti} to investigate the attitude tracking problem of rigid spacecraft systems.
In practical situation, the disturbance is inevitable during attitude maneuver of spacecraft systems and thus, it is crucial to take the disturbance into account in the attitude control problem of spacecraft systems (cf., \cite{kang1995nonlinear,luo2005inverse,chen2009attitude,chen2014attitude,dong2022anti,meng2023second}). In particular, a nonlinear $H_{\infty}$ method is proposed in \cite{kang1995nonlinear} to achieve disturbance attenuation in the attitude tracking problem of rigid spacecraft systems. 
The similar problem is further investigated in \cite{luo2005inverse} by an inverse optimal adaptive control approach, which achieves any desired level of $L_2$ disturbance attenuation. 
Later, the internal model approach is developed in \cite{chen2009attitude} to solve the attitude tracking problem of rigid spacecraft systems
subject to external disturbances with unbounded energy. In particular, disturbance rejection is achieved in \cite{chen2009attitude} where 
the disturbance can be any combination of finitely many sinusoidal functions with known frequencies.
Later, the result is further extended in \cite{chen2014attitude} to tackle the case where the frequencies of the disturbances are unknown. 
 


Performing more complex tasks, spacecrafts become larger with more attachments such as solar cell wings while the structural weight is required to be as small as possible. Thus, the flexibility issues become prominent and many researchers has investigated the attitude control problem of flexible spacecraft systems
(cf., \cite{xiao2011adaptive,zhong2017attitude,liu2017robust,xiao2011fault,xiao2020attitude,zhu2020disturbance,bao2023attitude}). In particular, 
an adaptive sliding mode control law is proposed in \cite{xiao2011adaptive} to tackle the attitude tracking problem of flexible spacecraft systems under external disturbances and the neural network is used to approximate the uncertain inertia matrix.
Fault-tolerant attitude stabilization control of uncertain flexible spacecraft systems in the presence of external disturbances and uncertain inertia matrix is investigated in \cite{xiao2011fault}. 
Recently, the attitude tracking problem is studied in \cite{zhong2017attitude} for flexible spacecraft systems subject to external disturbances with unbounded energy. To be precise, the 
attitude tracking and disturbance rejection are achieved where 
the disturbances can be a combination of sinusoidal signals with unknown amplitudes, phase angles and known frequencies in \cite{zhong2017attitude}.
In practice, the mass properties of the spacecraft is often uncertain or may change due to uncertain mass distribution, and spacecraft reconfiguration. Moreover, the frequencies of the external disturbance may be not measurable. Therefore, 
the problem that how to achieve attitude tracking of flexible spacecraft systems in the presence of the uncertain inertia matrix and external disturbances with unknown frequencies remains open.





In this paper, we investigate the attitude tracking problem of flexible spacecraft systems with uncertain inertia matrix and unknown external disturbances. In particular, the external disturbance can be any multi-tone sinusoidal signal with unknown frequencies, amplitudes and phase angles. We aim to solve the open problem, that is, to achieve exact attitude tracking of uncertain flexible spacecraft systems with disturbance rejection and vibration suppression of the flexible attachments by a sufficiently smooth control law. Technically, the uncertain inertia matrix, the unknown disturbance and the unmeasurable vibration modal give rise to great challenges to the solution of the problem. 
To overcome these challenges, we develop a class of novel nonlinear internal models based approach. 
First, we make a linearly parameterized manipulation to the unknown inertia matrix and propose a class of nonlinear internal models, which is used to compensate the unknown disturbance and partial nominal system dynamics. Based on the internal model principle, the attitude tracking
problem of uncertain flexible spacecraft systems is converted into a regulation problem of an augmented system composed of
the spacecraft system and the nonlinear internal model. Then, the regulation problem of the augmented system is solved by designing an auxiliary stabilization dynamic system, a dynamic compensator and an adaptive control law.
Specifically, the auxiliary dynamic system is introduced for auxiliary stabilization, which can circumvent the difficulty caused by the unmeasurable
modal variable and help the vibration suppression of the flexible attachments. The dynamic compensator is designed to compensate the uncertainty in the closed-loop system introduced by the nonlinear internal model. A linearly parameterized transformation is developed to deal with the challenges caused by the cross-coupling of the system uncertainty, the unknown disturbance frequencies and the system state.    
Finally, by introducing a series of
coordinate and input transformations, we propose an adaptive dynamic control law to achieve regulation of the augmented
system and thus leading to the solution to the attitude tracking problem. In addition, we analyze the convergence issue of the
estimated parameter to its true value by the persistently exciting (PE) condition. Finally, the effectiveness of the developed
approach is verified by its application to the attitude manoeuvre of a flexible spacecraft system in the presence of external
disturbances.

The remainder of this paper is organized as follows. The attitude tracking problem of uncertain flexible spacecraft systems subject to external disturbances is formulated in Section \ref{spfap}.
Section \ref{secmrs} presents our main results, which includes the nonlinear internal model design, nonlinear regulation, adaptive control and convergence analysis of estimated parameters.
An example is given in Section \ref{sae} to illustrate the effectiveness of the proposed control approach. Finally, some conclusions are
made to end this paper in Section \ref{con}.

{\bf Notation.} For $x_i \in \mathbb{R}^{n_i}$, $i=1,\dots,m$, $col(x_1,\dots,x_m)=[x_1^T,\dots,x_m^T]^T$.
For $A_i \in \mathbb{R}^{n\times p}$, $i=1,\dots,m$, $\mbox{block diag}(A_1,\dots,A_m)$ denotes a block matrix with its $i$th diagonal entry being $A_i$ and all other entries being zero.
$\otimes$ denotes the Kronecker product. 
%
For $A \in \mathcal{R}^{m \times n}$ and $B \in \mathcal{R}^{n \times p}$, $A \circ B$ denotes the Tracy-Singh product of $A$ and $B$ with
$(A \circ B)_{ij}= \sum_{k=1}^{n} (A_{ik}B_{kj})$, where $(A \circ B)_{ij}$  represents the (i, j)th element of $A \circ B$. 
For $a_i \in \mathbb{R}, i=1, \cdots, p$, diag$(a_1,\cdots,a_p)$ denotes the diagonal matrix with $i$th diagonal element being $a_i$ and all other elements being zero. For $C_i \in \mathbb{R}^{n \times m}$, $i=1,\cdots,p$, block diag$(C_1,\cdots,C_p)$ denotes a block diagonal matrix with its $i$th diagonal entry being $C_i$ and all other entries being zero. For $\varpi=[\varpi_1,\varpi_2,\varpi_3]^T \in \mathbb{R}^3$, define
$\varpi^{\times}=\left[
    \begin{array}{ccc}
        0 &-\varpi_3& \varpi_2\\
        \varpi_3 & 0 & -\varpi_1\\
        -\varpi_2& \varpi_1 & 0\\
    \end{array}
\right]$.



\section{\scshape Problem Formulation and Preliminaries}\label{spfap}

In this paper, we investigate motion control of the flexible spacecraft, whose model contains the attitude kinematics and the motion dynamics. The attitude  kinematic equation of the flexible spacecraft can be expressed as follows:
\begin{equation}\label{system1}
    \dot{q}=\frac{1}{2}
    \left[
        \begin{aligned}
           &q_v^{\times}+q_4I_3\\
           &-q_v^{T}
        \end{aligned}
    \right]\omega
\end{equation}
where $q=[q_v^T,q_4]^T$ represents the attitude quaternion of the flexible spacecraft with $q_v=[q_1,q_2,q_3]^T \in \mathbb{R}^3$ and $q_4 \in \mathbb{R}$, and $\omega=[\omega_1,\omega_2,\omega_3]^T \in \mathbb{R}^3$ is the angular velocity of the flexible spacecraft. The quaternion $q$ satisfies $q_v^T q_v+q_4^2=1$, which is called a unit quaternion. For a unit quaternion, the conjugate of $q$ is defined by $q^{-1}=[-q_v^T,q_4]^T$. The set of all unit quaternions is denoted by $\mathbb{Q}_u$. Compared with the Euler angles, the unit quaternion expression can avoid singularity. 


The motion equation of the flexible spacecraft can be expressed as follows \cite{monaco1985nonlinear}:
\begin{equation}\label{system2}
\begin{split}
J\dot{\omega}+\omega^{\times}J \omega+\delta \ddot{\eta}&=u+d\\
\ddot{\eta}+C \dot{\eta}+K \eta&=-\delta^T \dot{\omega}
\end{split}
\end{equation}
where $J \in \mathbb{R}^{3 \times3}$ is the symmetric positive definite inertia matrix and contains uncertainties,  $u \in \mathbb{R}^{3}$ is the control torque, and $\eta \in \mathbb{R}^n$ is the modal coordinate vector of the flexible attachments with $n$ being the number of the flexible modes. The variables $q, \omega, J, u$ and $\eta$ are all expressed in the body frame.  $\delta \in \mathbb{R}^{3 \times n}$ is the coupling matrix between the rigid body and flexible attachments,  $C \in \mathbb{R}^{n \times n}$ is the symmetric positive definite damping matrix and $K \in \mathbb{R}^{n \times n}$ is the symmetric positive definite stiffness matrix. The external disturbance $d(t)=[d_1(t),d_2(t),d_3(t)]^{\mathrm{T}}$ is assumed to be a class of persistent disturbances with unbounded energy which takes the following form:
\begin{equation}\label{raodong}
d_i(t)=C_{i0}+\sum_{j=1}^{N_i} C_{ij}{\rm sin}(b_{ij}t+l_{ij}), i=1,2,3
\end{equation}
where $C_{i0}$ are unknown step magnitudes, and $C_{ij}$,  
$b_{ij}$, and $l_{ij}$, $j=1,\cdots,N_i$ are unknown amplitudes, frequencies and phases, respectively.



The desired attitude is described by $q_{d}=\left[q_{d v}^{\mathrm{T}},q_{d 4}\right]^{\mathrm{T}}$  with $q_{d v}=\left[q_{d 1}, q_{d 2}, q_{d 3}\right]^{\mathrm{T}}$ and $q_{d}\in \mathbb{Q}_u$. In this paper, we consider a rest-to-rest attitude tracking problem, and thus the desired angular velocity $\omega_d=0$.
Now, the attitude tracking problem of the flexible spacecraft system subject to external disturbances can be described as follows.
\begin{Problem}\label{problem1}
    
    In the context of the uncertain and flexible spacecraft system consisting of \eqref{system1} and \eqref{system2}, with a desired attitude denoted as $q_d$, the objective is to devise a control law represented by
    \begin{equation}\label{controllaw}
        \begin{split}
            u=\phi(\chi,q,q_d,\omega,\omega_d)\\
            \dot{\chi}=\psi(\chi,q,q_d,\omega,\omega_d)
        \end{split}
    \end{equation}
where $\chi \in \mathbb{R}^{n_{\chi}}$ for some positive integer $n_{\chi}$, $\phi$ and $\psi$ are some sufficiently smooth functions, such that, for any disturbances $d(t)$ of the form \eqref{raodong} and any initial state $q(0)\in \mathbb{Q}_u$, $\omega(0) \in \mathbb{R}^{3}$, $\eta(0) \in \mathbb{R}^{n}$, and $ \chi(0) \in \mathbb{R}^{n_{\chi}}$, the solution of the closed-loop system composed of \eqref{system1}, \eqref{system2} and \eqref{controllaw} is bounded and satisfies:
\begin{equation}
\begin{split}
    &\lim _{t \rightarrow \infty}\left(q(t)-q_d(t)\right)=0\\
    &\lim _{t \rightarrow \infty}\left(\omega(t)-\omega_d(t)\right)=0\\
    &\lim _{t \rightarrow \infty}\eta(t)=0.
\end{split}
\end{equation}
\end{Problem}

As that in \cite{sidi1997spacecraft}, the quaternion error $q_{e}=\left[q_{e v}^{\mathrm{T}}, q_{e 4}\right]^{\mathrm{T}}$ with $q_{e v}=\left[q_{e 1}, q_{e 2}, q_{e 3}\right]^{\mathrm{T}}$ is defined as follows:
\begin{equation}\label{quer1}
\begin{split}
q_{e v}&=q_{d 4} q_{v}-q_{d v}^{\times} q_{v}-q_{4} q_{d v}\\
q_{e 4}&=q_{d v}^{\mathrm{T}} q_{d v}+q_{4} q_{d 4}.
\end{split}
\end{equation}

It is noted that the error $q_{e}$ is subject to the constraint $q_{e v}^{\mathrm{T}}q_{ev}+q_{e4}q_{e4}=1$. The angular velocity error $\omega_e=\omega$. Then, by \eqref{quer1}, the error dynamics can be put as follows:
\begin{subequations}\label{cls21}
\renewcommand\theequation{\theparentequation-\arabic{equation}}
\begin{equation}\label{dynamics}
 J \dot{\omega}_{e}=-\omega_{e}^{\times} J \omega_{e}+u+d-\delta \ddot{\eta}\end{equation}
\begin{equation}\label{modal}
 \ddot{\eta}+C \dot{\eta}+K \eta=-\delta^{\mathrm{T}} \dot{\omega}_{e}
\end{equation}
\begin{equation}\label{kinetmatics}
     \dot{q}_{e}=\frac{1}{2}\left[q_{e 4} I_{3}+q_{e v}^{\times} ;-q_{ev}^{\mathrm{T}}\right]\omega_{e}.
\end{equation}
\end{subequations}

It can be verified that the flexible spacecraft attitude tracks the desired attitude if and only if $q_{e v}=0$. Then, Problem \ref{problem1} is equivalent to the objective that design a smoothly robust control law of the form \eqref{controllaw} to achieve 
\begin{equation}\label{proconv1}
\lim _{t \rightarrow \infty}\left(q_{e v}(t), \omega_{e}(t), \eta(t)\right)=(0,0,0).
\end{equation}


\begin{Remark}
Problem \ref{problem1} is studied in \cite{zhong2017attitude} where the inertia matrix of the flexible spacecraft and the disturbance frequencies are known. In contrast, these variables are unknown in this paper, which causes great challenges to the solution of Problem \ref{problem1}.
In particular, the inertia matrix cannot be used for feedback control and the nonlinear internal model in \cite{zhong2017attitude} cannot be used directly to solve the current problem.
\end{Remark}

\section{Main Results}\label{secmrs}

\subsection{Nonlinear Internal Model Design}\label{sdo}

System \eqref{cls21} contains two parts of uncertainties, one is the parameter uncertainty of the inertia matrix $J$ and the other is the external disturbance $d(t)$. These two uncertain parts are coupled in the stability analysis, which gives rise to some challenges to the control law design.

To deal with these challenges, we first make a linearly parameterized manipulation to the inertia matrix $J$. In particular,
as in \cite{ahmed1998adaptive} and \cite{chen2009attitude}, for any vector $x=[x_1,x_2,x_3]^T \in \mathbb{R}^3$, there exists $Jx=L(x)[J_{11}\ J_{22}\ J_{33}\ J_{23}\ J_{13}\ J_{12}]^T$ where $L(x)=\left[
    \begin{array}{cccccc}
        x_1 &0 & 0& 0& x_3& x_2\\
        0 & x_2 &0 &x_3 &0 &x_1\\
        0& 0 & x_3 &x_2& x_1 & 0\\
    \end{array}
\right]$ and $J=\left[
    \begin{array}{ccc}
        J_{11} & J_{12} & J_{13}\\
        J_{12} & J_{22} &J_{23}\\
        J_{13}& J_{23} & J_{33} \\
    \end{array}
\right]$. When some elements of $J$ are unknown, there always exists an unknown vector $\mu \in \mathbb{R}^{n_{\mu}}$ with dimension $0 \leq n_{\mu} \leq 6$, and some known matrices $\bar{L}_{1} \in \mathbb{R}^{6 \times n_{\mu}}$ and $\bar{L}_{0} \in \mathbb{R}^{6 \times 1}$ such that $[J_{11}\ J_{22}\ J_{33}\ J_{23}\ J_{13}\ J_{12}]^T=\bar{L}_1\mu +\bar{L}_0$. Therefore, 
\begin{equation}\label{lpm0-ex}
J x=L_{1}(x) \mu+L_{0}(x)
\end{equation}
where $L_{0}(x)=L(x) \bar{L}_{0}$ and $L_{1}(x)=L(x) \bar{L}_{1}$. It follows that
\begin{equation}\label{lpm1}
\begin{split}
-\omega^{\times}_eJ\omega_e&=F_1(\omega_e)\mu+F_0(\omega_e)\\
F_i(\omega_e)&=-\omega^{\times}_eL_i(\omega_e), i=0,1.
\end{split}
\end{equation}
Thus, the system \eqref{dynamics} becomes:
\begin{equation}\label{dynamic1}
J\dot{\omega}_e=F_1(\omega_e)\mu+F_0(\omega_e)+u+d-\delta\ddot{\eta}.
\end{equation}

Drawn from the principle of output regulation theory,
(cf.,\cite{huang2004nonlinear},\cite{chen2015stabilization},\cite{isidori1990output}), we design a nonlinear internal model to tackle the unknown external disturbance $d(t)$. 
To do this, we use the notation $\sigma \in \mathbb{R}^{n_{\sigma}}$ to denote the vector consisting of all unknown frequencies of the external disturbance.
Then, for $i=1,2,3$, there exist positive integers $r_i$ and real numbers
 $a_{i1}(\sigma),\cdots,a_{ir_i}(\sigma)$ such that
\begin{equation*}
\begin{split}
\frac{d^{r_i}  d_{i}(t)}{d t^{r_i} }\!\!=\!\!a_{i 1}(\sigma) d_{i}(t)\!+\!a_{i 2}(\sigma) \frac{d d_{i}(t)}{d t}\!+\!\cdots \!+\!a_{i r_{i}}(\sigma) \frac{d^{\left(r_{i}-1\right)} d_{i}(t)}{d t^{\left(r_{i}-1\right)}}
\end{split}
\end{equation*}
where $a_{i1}(\sigma),a_{i2}(\sigma), \ldots, a_{ir_{i}}(\sigma)$ are some constants such that all roots of the polynomial $P_i(\lambda)=\lambda^{r_i}-a_{i1}(\sigma)-a_{i2}(\sigma)\lambda- \cdots-a_{ir_i}(\sigma) \lambda^{{r_i}-1}$ are distinct with zero real parts.

Let
\begin{equation}\label{matrixPhi}
\setlength{\arraycolsep}{0.5pt}
\Phi_{i}(\sigma)\!\!=\!\!\left[\begin{array}{ccccc}
0 & 1 & 0 & \cdots & 0 \\
0 & 0 & 1 & \cdots & 0 \\
\vdots & \vdots & \vdots & \cdots & \vdots \\
0 & 0 & 0 & \cdots & 1 \\
a_{i, 1}(\sigma) & a_{i, 2}(\sigma) & a_{i, 3}(\sigma) & \cdots & a_{i, r_{i}}(\sigma)
\end{array}\right], 
\Psi_{i}\!\!=\!\!\left[\begin{array}{c}
1 \\
0 \\
\vdots \\
0
\end{array}\right]^{\mathrm{T}}
\end{equation}
and $T_i$ be any nonsingular matrix of dimension $r_i$, $i=1,2,3$. 

Defining
\begin{equation}
\varrho_{i}(t)=\left[\begin{array}{llll}
d_{i}(t) & \dot{d}_{i}(t) & \cdots & \left.\frac{d^{\left(r_{i}-1\right)} d_{i}(t)}{d t^{\left(r_{i}-1\right)}}\right]^{\mathrm{T}}
\end{array}\right.
\end{equation}
we can show that
\begin{equation}\label{varrho_idefintion}
\begin{split}
\dot{\varrho}_{i}(t)&=\Phi_{i}(\sigma)\varrho_i(t)\\
\quad d_{i}(t)&=\Psi_{i}  \varrho_{i}(t).
\end{split}
\end{equation}

Choose $M_{i} \in \mathbb{R}^{r_{i} \times r_{i}}$ and $N_{i} \in \mathbb{R}^{r_{i}}$ as a pair of controllable matrices where $M_{i}$ is Hurwitz. Furthermore, since $(M_i,N_i)$ is controllable and $(\Phi_i(\sigma),\Psi_i)$ is observable, there exists a nonsingular matrix $T_{i}(\sigma)$ satisfying the Sylvester equation \cite{nikiforov1998adaptive}:
\begin{equation}\label{Mdefinition}
T_{i}(\sigma) \Phi_{i}(\sigma)-M_{i} T_{i}(\sigma)=N_{i} \Psi_{i}.
\end{equation}
Let $\theta_i=-T_i(\sigma) \varrho_i$. Then,
\begin{equation}\label{steadss}
\begin{split}
    \dot{\theta}_{i}(t)&=T_{i}(\sigma) \Phi_{i}(\sigma)  T_{i}^{-1}(\sigma) \theta_{i}(t)\\
    \quad d_{i}(t)&=-\Psi_{i} T_{i}^{-1}(\sigma) \theta_{i}(t).
\end{split}
\end{equation}

Define
$\theta =\operatorname{col}\left(\theta_{1}, \theta_{2}, \theta_{3}\right),
T(\sigma)=\text{block}\operatorname{diag}(\mathrm{T}_{1}(\sigma), \mathrm{T}_{2}(\sigma),$ $ \mathrm{T}_{3}(\sigma)),
\Phi(\sigma)\!\!=\!\!\text{block}\!\operatorname{diag}\left(\mathrm{\Phi}_{1}(\sigma), \mathrm{\Phi}_{2}(\sigma), \mathrm{\Phi}_{3}(\sigma)\right),
\mathrm{X}\!\!=\!\!\text { block }\! \operatorname{diag}(\mathrm{X}_{1}, $ $\mathrm{X}_{2}, \mathrm{X}_{3}), \quad \mathrm{X}=M, N,\Psi .$
Then, we can design a nonlinear internal model as follows:
\begin{equation}\label{internal}
\dot{v}=Mv+Nu+NF_0(\omega_e)-MNL_0(\omega_e)
\end{equation}
where $v \in \mathbb{R}^r$ with $r=r_1+r_2+r_3$. 


Introducing the following coordinate transformation:
\begin{equation}\label{v_and_checkv}
    \begin{split}
        \check{v}&=v-\theta-N\delta\dot{\eta}-NJ\omega_e\\
        \tilde{u}&=u+F_0(\omega_e)
    \end{split}
\end{equation}
%
we obtain an augmented system as follows:
\begin{equation}\label{augsys1}
\begin{split}
J\dot{\omega}_e=&F_1(\omega_e)\mu+\tilde{u}+\Psi T^{-1}(\sigma )(\check{v}+N\delta \dot{\eta}+NJ\omega_e)\\
&-\Psi T^{-1}(\sigma )v-\delta\ddot{\eta}\\
\dot{\check{v}}=&M\check{v}+MN\delta \dot{\eta}+(MNL_1(\omega_e)-NF_1(\omega_e))\mu\\
\ddot{\eta}=&-C \dot{\eta}-K \eta-\delta^{\mathrm{T}} \dot{\omega}_{e}\\
\dot{q}_{e}=&\frac{1}{2}\left[q_{e 4} I_{3}+q_{e v}^{\times} ;-q_{ev}^{\mathrm{T}}\right]\omega_{e}.
\end{split}
\end{equation}




Now, the attitude tracking problem of system \eqref{cls21} can be converted into a regulation problem of the system \eqref{augsys1}.
The result is summarized as follows.
\begin{Lemma}\label{Lammaaux1}
By the transformation \eqref{v_and_checkv} and the closed system \eqref{augsys1}, suppose that there exists a dynamic control law of the form 
    \begin{equation}\label{controllaw1}
        \begin{split}
            \tilde{u}=h(\bar{\chi},q_e,\omega_e)\\
            \dot{\bar{\chi}}=k(\bar{\chi},q_e,\omega_e)
        \end{split}
    \end{equation}
where $\bar{\chi} \in \mathbb{R}^{n_{\bar{\chi}}}$ is the compensator state, $h$ and $k$ are some sufficiently smooth functions, such that for any initial state $q(0)\in \mathbb{Q}_u$, $\omega(0)\in \mathbb{R}^3$, $\eta(0)\in \mathbb{R}^n$, and $\bar{\chi}(0)\in \mathbb{R}^{n_{\bar{\chi}}}$, the state of the system composed of
\eqref{augsys1} and \eqref{controllaw1}
is bounded and $\lim _{t \rightarrow \infty}\left(q_{e v}(t), \omega_{e}(t), \eta(t)\right)=(0,0,0)$. Then, the attitude tracking Problem \ref{problem1} is solved by the control law of the form 
\begin{equation}
\begin{split}
 u&=h(\bar{\chi},q_e,\omega_e)-F_0(\omega_e)\\
 \dot{v}&=Mv+Nu+NF_0(\omega_e)-MNL_0(\omega_e)\\
 \dot{\bar{\chi}}&=k(\bar{\chi},q_e,\omega_e).
\end{split}
\end{equation}
\end{Lemma}



\begin{Remark}
The dynamic compensator \eqref{internal} is inspired by the internal model design approach in \cite{huang2004nonlinear}, \cite{chen2009attitude} and \cite{zhong2017attitude}. However, the existing internal models in these works cannot be used to solve the current problem. In particular, when the inertia matrix is known, a nonlinear internal model of the form $\dot{v}=f(v,\omega_e, J)$ is proposed in \cite{zhong2017attitude}. It is noted that the inertia matrix $J$ is uncertain in the work and thus the internal model in \cite{zhong2017attitude} cannot be applied. Therefore, we design a novel nonlinear internal model which is different from the existing ones. 
\end{Remark}

\par

\subsection{Nonlinear Regulation}\label{smr}

In this subsection, we develop a dynamic control law of the form \eqref{controllaw1} to solve the nonlinear regulation problem of system \eqref{augsys1} in Lemma \ref{Lammaaux1}.
To achieve the nonlinear regulation is rather challenging. First, the mechanical actuator of the flexible spacecraft is underactuated and thus the modal variable is 
not accessible for feedback control. It makes the vibration suppression of the flexible attachments difficult. Second, 
the inertia matrix is allowed to contain unknown parameters and thus cannot be used in the control law, which fails the feedfoward design in the control law. Third, the uncertainty caused by the inertia matrix and that by the unknown disturbance frequencies are cross-coupled with the system state in the system dynamics which precludes the regulation of the system.

\subsubsection{Auxiliary Stabilization}\label{sssas}

To achieve  vibration suppression without modal variable measurement, as in \cite{zhong2017attitude}, we introduce an auxiliary system as follow:
\begin{equation}\label{etas1}
\dot{z}=Az+H(\omega_e)
\end{equation}
where 
\begin{equation}\label{etas1matr}
\begin{split}
z&=\left[\begin{array}{l}
z_{1} \\
z_{2}
\end{array}\right], z_2=\eta \\
A&=\left[\begin{array}{cc}
0 & I_{n} \\
-K & -C
\end{array}\right],
H\left(\omega_{e}\right)=\left[\begin{array}{l}
0 \\
-\delta^{\mathrm{T}} \omega_{\mathrm{e}}
\end{array}\right].
\end{split}
\end{equation}



Substitute \eqref{etas1} into \eqref{augsys1} gives
\begin{equation}
\begin{aligned}
J\dot{\omega}_e=&F_1(\omega_e)\mu+\tilde{u}+\Psi T^{-1}(\sigma )(\check{v}+N\delta \dot{\eta}+NJ\omega_e)\\&-\Psi T^{-1}(\sigma )v-\delta\ddot{z}_2\\
=&F_1(\omega_e)\mu+\tilde{u}+\Psi T^{-1}(\sigma)\check{v}+\Psi T^{-1}(\sigma)N(L_1(\omega_e)\mu\\
&+L_0(\omega_e))+\Psi T^{-1}(\sigma)N\delta(-Kz_1-Cz_2-\delta^T \omega_e)\\
&-\Psi T^{-1}(\sigma )v-\delta CKz_1-\delta(C^2-K)z_2-\delta C \delta^T\omega_e\\
&+\delta \delta^{T} \dot{\omega}_e.
\end{aligned}
\end{equation}

Furthermore, defining $J_{mb}=(J-\delta \delta^{T})$, we have 
\begin{equation}\label{dynamicsJmb}
\begin{split}
J_{mb}\dot{\omega}_e=
&F_1(\omega_e)\mu+\tilde{u}+\Psi T^{-1}(\sigma)\check{v}\\
&+\Psi T^{-1}(\sigma)N(L_1(\omega_e)\mu+L_0(\omega_e))\\
&-\Psi T^{-1}(\sigma)N\delta\delta^T \omega_e-\Psi T^{-1}(\sigma )v\\
&-(\delta(C^2-K)+\Psi T^{-1}(\sigma)N\delta C)z_2 \\
&-(\delta CK+\Psi T^{-1}(\sigma)N\delta K)z_1\\
&-\delta C \delta^T\omega_e.
\end{split}
\end{equation}

\subsubsection{Dynamic Compensation}\label{sssas}

To deal with the uncertainty by the inertia matrix in nonlinear regulation, we  design a dynamic compensator as follows: 
\begin{equation}\label{dyncomp1}
\dot{\zeta}=M\zeta+(MNL_1(\omega_e)-NF_1(\omega_e))
\end{equation}
where $\zeta \in \mathbb{R}^{r\times n_{\mu} }$.

Define the coordinate transformation as follows: 
\begin{equation}\label{inertcord1}
\hat{v}=\check{v}-\zeta \mu.
\end{equation}
Then, by \eqref{augsys1}, \eqref{etas1}  and \eqref{dynamicsJmb}, we have
\begin{equation}\label{etas2}
\begin{split}
\dot{\hat{v}}=&\dot{\check{v}}-\dot{\zeta}\mu\\
=&M\check{v}+(MNL_1(\omega_e)-NF_1(\omega_e))\mu
+MN\delta\dot{\eta}\\
&-(M\zeta+(MNL_1(\omega_e)-NF_1(\omega_e)))\mu\\
=&M\hat{v}+MN\delta\dot{\eta}\\
=&M\hat{v}-MN\delta Cz_2-MN\delta Kz_1-MN\delta \delta^{T}\omega_e\\
\end{split}
\end{equation}
and 
\begin{equation}\label{syseqco1}
\begin{split}
J_{mb}\dot{\omega}_e=
&F_1(\omega_e)\mu+\tilde{u}+\Psi T^{-1}(\sigma)\hat{v}+\Psi T^{-1}(\sigma)\zeta\mu\\
&+\Psi T^{-1}(\sigma)N(L_1(\omega_e)\mu+L_0(\omega_e))\\
&-\Psi T^{-1}(\sigma)N\delta\delta^T \omega_e-\Psi T^{-1}(\sigma )v\\
&-(\delta(C^2-K)+\Psi T^{-1}(\sigma)N\delta C)z_2\\
&-(\delta CK+\Psi T^{-1}(\sigma)N\delta K)z_1-\delta C \delta^T\omega_e.
\end{split}
\end{equation}

\begin{Remark}
By  \eqref{etas2}, it can be found that  the uncertainty term $(MNL_1(\omega_e)-NF_1(\omega_e))\mu$ in \eqref{augsys1} is eliminated by introducing the dynamic compensator system $\zeta$ and the coordinate transformation \eqref{inertcord1}. That is, the uncertainty parameter caused the inertia matrix and that by the external disturbance in the subsystem $\check{v}$ are decoupled, which makes it possible to achieve nonlinear regulation.
\end{Remark}


\subsubsection{Adaptive Regulation}\label{sssas}

It can be found that the uncertainty caused by the inertia matrix and that by the unknown disturbance frequencies in \eqref{syseqco1} are cross-coupled with the state variable, which results in great challenges for the control law design.
To deal with the uncertainty,
we first make a linearly parameterized manipulation as follows:
\begin{equation}\label{Edefinitionlabel}
\begin{aligned}
& \Psi T^{-1}(\sigma)=E_0+\sum_{j=1}^{\ell} E^j \Omega^j(\sigma)=E_0+E\left(\Omega^{\sigma} \otimes I_r\right) \\
& E:=\left[E^1, \cdots, E^{\ell}\right], \Omega^{\sigma}:=\left[\Omega_1^{\sigma}, \cdots, \Omega_{\ell}^{\sigma}\right]^{\top}
\end{aligned}
\end{equation}
where $\ell$ is an integer satisfies $\ell \geq 1$, $E_0 \in \mathbb{R}^{3 \times r}$ and $E^{j} \in \mathbb{R}^{3 \times r}$ are constant matrices, and $\Omega_j^{\sigma} \in \mathbb{R}$ is a smooth function of the unknown parameter $\sigma$, where $j=1,\cdots,\ell$.

Furthermore, the cross-coupling term can be put as follows:
\begin{equation} 
\begin{aligned}
\Psi T^{-1}\left(\sigma\right) \zeta \mu  =&(E_0+E\left(\Omega^{\sigma} \otimes I_r\right)) \zeta \mu \\
 =&E_0 \zeta \mu+(E \circ \zeta)(\Omega^{\sigma} \otimes \mu)
\end{aligned}
\end{equation}
\begin{equation}
\begin{aligned}
&\Psi T^{-1}(\sigma)\left(N L_0(\omega_e)-v\right)\\
&=(E_0+E\left(\Omega^{\sigma} \otimes I_r\right))\left(N L_o(\omega_e)-v\right)\\
&=E_0\left(N L_o(\omega_e)-v\right) +E \circ\left(N L_o(\omega_e)-v\right) \Omega^{\sigma} \\
\end{aligned}
\end{equation}
and
\begin{equation}
\begin{aligned}
&\Psi T^{-1}(\sigma) N L_1(\omega_e) \mu  \\
&=(E_0+E\left(\Omega^{\sigma} \otimes I_r\right)) N L_1(\omega_e) \mu\\
&=E_0N L_1(\omega_e) \mu+E \circ\left(N L_1(\omega_e)\right)(\Omega^{\sigma} \otimes \mu) .
\end{aligned}
\end{equation}
Define 
\begin{equation}\label{matrixrhoR}
\begin{split}
\rho(\omega_e,\zeta,v)=&[\rho_1 \quad \rho_2 \quad \rho_3]
\end{split}
\end{equation}
where $\rho_1=F_1(\omega_e)+E_0(NL_1(\omega_e)+\zeta)$, $\rho_2=E \circ(\zeta+NL_1(\omega_e))$ and $\rho_3=E \circ (NL_0(\omega_e)-v)$. 
Denote
\begin{equation}
\begin{split}
R(\sigma,\mu)=[\mu,\Omega^{\sigma} \otimes \mu, \Omega^{\sigma}]^{\mathrm{T}}.
\end{split}
\end{equation}
Then, introducing the input transformation as follows:
\begin{equation}
\tilde{u}=\check{u}-E_0N L_o(\omega_e)+E_0v+\delta C \delta^{T} \omega_e
\end{equation}
we have
\begin{equation}\label{close23}
\begin{split}
\dot{\hat{v}}=&M\hat{v}-MN\delta Cz_2-MN\delta Kz_1-MN\delta \delta^{T}\omega_e\\
\dot{\omega}_e=&J_{mb}^{-1}(\Psi T(\sigma)^{-1}\hat{v}+\rho(\omega_e,\zeta,v)R(\sigma,\mu)\\
&+\check{u}-(\delta(C^2-K)+\Psi T^{-1}(\sigma)N\delta C)z_2\\
&-(\delta CK+\Psi T^{-1}(\sigma)N\delta K)z_1-\Psi T^{-1}(\sigma)N \delta \delta^T \omega_e)\\
\dot{q}_{e}=&\frac{1}{2}\left[q_{e 4} I_{3}+q_{e v}^{\times} ;-q_{e v}^{\mathrm{T}}\right] \omega_{e}\\
\dot{z}=&Az+G(\omega_e).
\end{split}
\end{equation}

We design an updated adaptive control law as follows:
\begin{equation}\label{adplaw1}
\dot{\hat{R}}=k\rho^T(\omega_e,\zeta,v) \omega_e
\end{equation}
where the state $\hat{R}$ is used to estimate $R(\sigma,\mu)$, and $k$ is any positive real number.

Based on \eqref{adplaw1}, we propose the virtual control law $\check{u}$ as follows:
\begin{equation}\label{uequ}
\check{u}=-k_1q_{ev}-k_2\omega_e-\rho(\omega_e,\zeta,v)\hat{R}
\end{equation}
where $k_1$ and $k_2$ are some positive real numbers to be determined. 

Now, we can show that the adaptive control law composed of \eqref{adplaw1} and \eqref{uequ} can solve the regulation problem of system \eqref{close23}, which leads to the solution to Problem \ref{problem1}.

\begin{Theorem}\label{them1}
The attitude tracking problem of the uncertain flexible spacecraft system composed of \eqref{system1} and \eqref{system2}
can be solved by the following adaptive control law:
\begin{equation}\label{ctrf1}
\begin{split}
        u=&-k_1q_{ev}-k_2\omega_e-\rho(\omega_e,\zeta,v) \hat{R}-E_{0}NL_0(\omega_e)\\
        &+\delta C \delta^{T} \omega_e+E_{0}v-F_0(\omega_e)\\
        \dot{v}=&Mv+Nu+NF_0(\omega_e)-MNL_0(\omega_e)\\
        \dot{\hat{R}}=&k\rho(\omega_e,\zeta,v)^T \omega_e\\
        \dot{\zeta}=&M \zeta+\left(M N L_{1}\left(\omega_{e}\right)-N F_{1}\left(\omega_{e}\right)\right).
\end{split}
\end{equation}
\end{Theorem}

\begin{Proof}
Since $C$ and $K$ are symmetric positive definite matrices, $A$ in \eqref{etas1matr} is always Hurwitz. Thus, for any $p>0$, 
there always exists positive definite matrix $P=P^T \in \mathbb{R}^{2 n \times 2 n}$ such that
\begin{equation}\label{Pdefin}
P A+A^{\mathrm{T}} P=-p I.
\end{equation}
Since $P$ is positive definite, there exist positive real numbers $\beta_1$ and $\beta_2$ such that $ \beta_{1}  x^T  x \leq x^T P x \leq  \beta_{2}  x^T  x $ for any $x\in \mathbb{R}^{2n}$.

Since $M$ is Hurwitz, for any $s>0$, there exists a positive definite matrix $S=S^T \in \mathbb{R}^{r \times r}$ such that 
\begin{equation}\label{Sdefin}
SM+M^TS=-sI_r.
\end{equation}

Define $\tilde{R}=\hat{R}-R$ and 
$V(t)=V_1(t)+V_2(t)+V_3(t)$, 
where 
\begin{equation}
\begin{split}
V_1(t)&=z^T Pz\\
V_2(t)&=k_1((q_{e4}-1)^2+q_{ev}^{T}q_{ev})+\frac{1}{2}\omega_e^{T}J_{mb}\omega_e +\frac{1}{2k}\tilde{R}^T\tilde{R}\\
V_3(t)&=\frac{1}{2}\hat{v}^T S \hat{v}.
\end{split}
\end{equation}

First, by \eqref{Pdefin}, the derivative of $V_1$ along the trajectory of $z$ subsystem of \eqref{close23}  satisfies
\begin{equation}
\begin{aligned}\label{v1}
\dot{V}_{1}(t)&=\dot{z}^{\mathrm{T}} P z+z^{\mathrm{T}}P\dot{z}\\
&=-p\|z\|^{2}+2 z^{\mathrm{T}} P H\left(\omega_{e}\right)\\ 
&\leq-\left(p -\frac{1}{\epsilon_{1}}\right)\|z\|^{2}+\epsilon_{1}  \beta_{2}^{2}\|\delta\|^{2}\left\|\omega_{e}\right\|^{2}
\end{aligned}
\end{equation}
where $\epsilon_1 $ is some positive real number, and  $\left\|H\left(\omega_{e}\right)\right\| \leq\|\delta\|\left\|\omega_{e}\right\|$.

Second, by \eqref{uequ}, the derivative of $V_2$ along  the trajectory of $\mbox{col}(q_e,\omega_e)$ subsystem of the system \eqref{close23} and $\hat{R}$ system \eqref{adplaw1} satisfies
\begin{equation}\label{v2}
\begin{split}
\dot{V_2}(t)=&k_1\omega_e^{T}q_{ev}+\omega_e^{T}(\Psi T^{-1}(\sigma)\hat{v}+\rho(\omega_e,\zeta,v)R(\sigma,\mu)\\
&+\check{u}-(\delta(C^2-K)+\Psi T^{-1}(\sigma)N\delta C)z_2\\
&-(\delta CK+\Psi T^{-1}(\sigma)N\delta K)z_1-\Psi T^{-1}(\sigma)N \delta \delta^T \omega_e)\\
&+\tilde{R}^T\rho^{T}(\omega_e,\zeta,v)\omega_e\\
=&-k_2 \omega_e^T \omega_e-\omega_e^T\alpha\omega_e-\omega_e^T \alpha_1 z_2-\omega_e^T \alpha_2 z_1+\omega_e^T \alpha_3 \hat{v}\\
 \leq&  -(k_2\!-\! \frac{1}{\epsilon_2}) ||\omega_e||^2\!+\!\frac{1}{4}||\omega_e^T \alpha||^2 \epsilon_2 \!+\!\frac{1}{4}||\omega_e^{T} \alpha_1||^{2}\epsilon_3+\frac{||z_2||^2}{\epsilon_3}
\\
 &+\frac{1}{4} ||\omega_e^{T} \alpha_2||^{2} \epsilon_4+\frac{||z_1||^2}{\epsilon_4}
+\frac{1}{4} ||\omega_e^{T} \alpha_3||^{2} \epsilon_5+\frac{||\hat{v}||^2}{\epsilon_5}
\end{split}
\end{equation}
where $\epsilon_1$, $\epsilon_2$, $\epsilon_3$, $\epsilon_4$, $\epsilon_5$ are positive real numbers, $\alpha=\Psi T^{-1}(\sigma)N \delta \delta^T$, $\alpha_1=[\delta(C^2-K)+\Psi T^{-1}(\sigma)N\delta C]$, $\alpha_2=[\delta CK+\Psi T^{-1}(\sigma)N\delta K]$, and $\alpha_3=\Psi T(\sigma)^{-1}$.

Third, noting \eqref{Sdefin}, we have that the derivative of $V_3$ along the trajectory of $\hat{v}$ subsystem of \eqref{close23}  satisfies
\begin{equation}\label{v3}
\begin{split}
\dot{V_3}(t)=&\frac{1}{2}\hat{v}^T S \dot{\hat{v}}+\frac{1}{2}\dot{\hat{v}}^T S \hat{v}\\
=&\hat{v}^T S(M\hat{v}-MN\delta Cz_2-MN\delta Kz_1-MN\delta \delta^{T}\omega_e)\\
=&-\frac{1}{2}s||\hat{v}||^2-\hat{v}^T \alpha_4 z_1-\hat{v}^T \alpha_5 z_2-\hat{v}^T \alpha_6 \omega_e\\
 \leq& -(\frac{1}{2}s-\frac{1}{4}||\alpha_4||^2\epsilon_7-\frac{1}{4}||\alpha_5||^2\epsilon_6-\frac{1}{4}||\alpha_6||^2\epsilon_8)||\hat{v}||^2\\
&+\frac{1}{\epsilon_7}||z_1||^2+\frac{1}{\epsilon_6}||z_2||^2+\frac{1}{\epsilon_8}||\omega_e||^2
\end{split}
\end{equation}
where $\alpha_4=SMN\delta K$, $\alpha_5=SMN\delta C$, and $\alpha_6=SMN\delta \delta^{T}$.

By \eqref{v1}, \eqref{v2} and \eqref{v3}, we have
\begin{equation}
\begin{split}
\dot{V}(t)
\leq &-(p-\frac{1}{\epsilon_1}-\frac{1}{\epsilon_4}-\frac{1}{\epsilon_7})||z_1||^2-(k_2-\frac{1}{\epsilon_2}-\frac{1}{4}||\alpha||^2 \epsilon_2\\
&-\frac{1}{4}||\alpha_1||^2\epsilon_3-\frac{1}{4}||\alpha_2||^2\epsilon_4-\frac{1}{4}||\alpha_3||^2\epsilon_5-\frac{1}{\epsilon_8}\\
&-\epsilon_1 \beta_2^2||\delta||^{2})||\omega_e||^2 -(\frac{1}{2}s-\frac{1}{\epsilon_5}-\frac{1}{4}||\alpha_4||^2 \epsilon_7\\
&-\frac{1}{4}||\alpha_5||^2 \epsilon_6 -\frac{1}{4}||\alpha_6||^2 \epsilon_8)||\hat{v}||^2\\
&-(p-\frac{1}{\epsilon_1}-\frac{1}{\epsilon_3} -\frac{1}{\epsilon_6})||z_2||^2.\\
\end{split}
\end{equation}

Choosing
\begin{equation}\label{parameterend12}
\begin{split}
s\geq&  \frac{2}{\epsilon_5}+\frac{1}{2}||\alpha_4||^2 \epsilon_7+\frac{1}{2}||\alpha_5||^2 \epsilon_6 +\frac{1}{2}||\alpha_6||^2 \epsilon_8+1\\
%
p\geq& \max\{\frac{1}{\epsilon_1}+\frac{1}{\epsilon_3}+\frac{1}{\epsilon_6}+1,\frac{1}{\epsilon_1}+\frac{1}{\epsilon_4}+\frac{1}{\epsilon_7}+1\}\\
k_1\geq& 1\\
k_2\geq& \frac{1}{\epsilon_2}+\frac{1}{4}||\alpha||^2 \epsilon_2+\frac{1}{4}||\alpha_1||^2\epsilon_3+\frac{1}{4}||\alpha_2||^2\epsilon_4\\
&+\frac{1}{4}||\alpha_3||^2\epsilon_5+\frac{1}{\epsilon_8}+\epsilon_1 \beta_2^2||\delta||^{2}+1 \\
\end{split}
\end{equation}
leads to  
\begin{equation}
\begin{split}
\dot{V}(t)
&\leq- ||z_1||^2- ||\omega_e||^2-||\hat{v}||^2-||z_2||^2.\\
\end{split}
\end{equation}
Since $V(t) \geq 0$, $z(t)$, $\omega_e(t)$, $q_{e}(t)$, $\tilde{R}(t)$, and $\hat{v}(t)$ are all bounded. Furthermore, by \eqref{kinetmatics}, \eqref{etas2}, and \eqref{etas1}, we have that $\dot{q}_{ev}(t)$, $\dot{\hat{v}}(t)$ and $\dot{z}(t)$ are all bounded. 

By \eqref{close23} and  \eqref{uequ}, we have
\begin{equation}\label{closedynamicu}
    \begin{split}
       \dot{\omega}_e=&J_{mb}^{-1}(\Psi T(\sigma)^{-1}\hat{v}-\rho(\omega_e,\zeta,v)\tilde{R}(\sigma,\mu)-k_1q_{ev}\\
       &-k_2\omega_e-(\delta(C^2-K)+\Psi T^{-1}(\sigma)N\delta C)z_2\\
       &-(\delta CK+\Psi T^{-1}(\sigma)N\delta K)z_1-\Psi T^{-1}(\sigma)N \delta \delta^T \omega_e)
    \end{split}
\end{equation}
which implies the boundedness of $\dot{\omega}_e(t)$.
Thus, we have that $\ddot{V}(t)$ is bounded. By Barbalat's lemma \cite{slotine1991applied}, we have 
$\lim_{t \to \infty} \dot{V}(t)=0$, which leads to
$$\lim_{t \to \infty} (z(t),\omega_e(t),\hat{v}(t))=(0,0,0).$$
It implies that 
$\lim_{t \to \infty} \eta(t)=0$.

By \eqref{closedynamicu}, $\ddot{\omega}_e(t)$ is bounded. Again, by Barbalat's lemma, $\lim_{t \to \infty} \dot{\omega}_e(t)=0$, 
which together with \eqref{closedynamicu} implies that $\lim_{t \to \infty}q_{ev}(t)=0$. 
The proof is thus completed by invoking Lemma \ref{Lammaaux1}.
\end{Proof}

\subsection{Convergence of the Estimated Parameters}\label{cep}
In this section, we   analyze the convergence of the estimated parameters to its true value based on the persistently exciting condition.
For this purpose, we introduce the following persistently exciting property of a signal.
\begin{Definition}\label{definitionpe}
\cite{ioannou1996robust}
A time function $f(t):[0, \infty) \mapsto \mathbb{R}^n$ is said to be persistently exciting (PE) if there exist $t_o, T_o, \vartheta>0$, such that $\left(1 / T_0\right) \int_t^{t+T_0}\left|c^{\top} f(s)\right| d s \geqslant \vartheta, \forall t \geqslant t_0$ for any unit vector $c \in \mathbb{R}^n$. A time function matrix $f(t):[0, \infty) \mapsto \mathbb{R}^{n \times m}$ is said to be PE if there exists a vector $b \in \mathbb{R}^m$ such that the vector $f(t) b$ is PE.
\end{Definition}
For convenience of subsequent convergence analysis, we define the signals as follows:
\begin{equation}\label{coveq1}
    \begin{split}
        A(t)&=e^{Mt}A_0,
        B(t)=E_0A(t)\\
    \end{split}
\end{equation}
where $M$ is defined in \eqref{Mdefinition}, $E_0$ is defined in \eqref{Edefinitionlabel}, and $A_0 \in \mathbb{R}^{r \times n_{\mu}}$ is an arbitrary real matrix.

The convergence of the estimated parameters is established by the theorem as follows.
\begin{Theorem}\label{theoremofpe}
    Let 
    \begin{equation}\label{equationy12}
        y(t):=\left[\begin{array}{lll}B(t) & E \circ A(t) & E \circ\left(T(\sigma) \varrho(t)-A(t) \mu\right)\end{array}\right]
    \end{equation}
    where $E$ is defined in \eqref{Edefinitionlabel}, and $A(t)$ and $B(t)$ are defined in \eqref{coveq1}.
    Then, $\lim _{t \rightarrow \infty} y(t)(\hat{R}(t)-R(\sigma, \mu))=0$. Moreover, if $y(t)$ is PE, then, $\lim _{t \rightarrow \infty}(\hat{R}(t)-R(\sigma, \mu))=0$.
\end{Theorem}
\begin{Proof}
As in  Theorem \ref{them1}, 
by \eqref{close23} and \eqref{uequ},
\begin{equation}\label{dotomegaequation}
    \begin{split}
       \dot{\omega}_e=&J_{mb}^{-1}(\Psi T(\sigma)^{-1}\hat{v}-\rho(\omega_e,\zeta,v)\tilde{R}(\sigma,\mu)-k_1q_{ev}\\
       &-k_2\omega_e-(\delta(C^2-K)+\Psi T^{-1}(\sigma)N\delta C)z_2\\
       &-(\delta CK+\Psi T^{-1}(\sigma)N\delta K)z_1-\Psi T^{-1}(\sigma)N \delta \delta^T \omega_e).\\ 
    \end{split}
\end{equation}

Note that by Theorem \ref{them1}, $\lim_{t \to \infty} (\dot{\omega}_e(t), \hat{v}(t), \omega_e(t), $\\ $z(t), q_{ev}(t))=0$, which together with \eqref{dotomegaequation} implies that 
\begin{equation}\label{pecondition0}
\lim_{t \to \infty}\rho(\omega_e(t),\zeta(t),v(t))\tilde{R}(t)=0.
\end{equation}

Next, we show that   
\begin{equation}\label{pecondition}
\lim_{t \to \infty} (y(t)-\rho (\omega_e(t),\zeta(t),v(t)))=0.
\end{equation}

By \eqref{coveq1} and \eqref{dyncomp1}, we have 
\begin{equation}
    \begin{split}
       \dot{A}(t)-\dot{\zeta }(t)&=MA(t)-M\zeta(t)-(MNL_1(\omega_e)-NF_1(\omega_e))\\
       &=M(A(t)-\zeta(t))-(MNL_1(\omega_e)-NF_1(\omega_e))\\
    \end{split}
\end{equation}
where $\zeta(t)$ is defined in \eqref{dyncomp1}. Since $\lim_{t \to \infty} \omega_e(t)=0$ and $M$ is Hurwitz matrix, we have
\begin{equation}\label{peimde1}
\lim_{t \to \infty} (A(t)-\zeta(t))=0.
\end{equation}

By  \eqref{v_and_checkv} and \eqref{inertcord1},  
\begin{equation}\label{sdf15}
    \begin{split}
    v&=\check{v}+\theta+\delta\dot{\eta}+NJ\omega_e.
    \end{split}
\end{equation}

By Theorem \ref{them1},  $\lim_{t \to \infty} \hat{v}(t)=\lim_{t \to \infty}(\check{v}(t)-\zeta(t) \mu)=0$, which together with \eqref{peimde1} and \eqref{sdf15} gives 
\begin{equation}\label{rho3y3}
    \begin{split}
        &\lim_{t \to \infty} (T(\sigma) \varrho(t)-A(t) \mu+v(t))\\
        =& \lim_{t \to \infty} (T(\sigma) \varrho(t) -A(t)\mu +\check{v}(t)+\theta(t)+NJ\omega_e(t)+N\delta\dot{\eta}(t))\\
        =&\lim_{t \to \infty}(\check{v}(t)-A(t)\mu)\\
        =&\lim_{t \to \infty}(\zeta(t) \mu-A(t) \mu)=0
    \end{split}
\end{equation}
where $\varrho=col(\varrho_1,\varrho_2,\varrho_3)$, and $\varrho_i$ is defined in \eqref{varrho_idefintion}, $i=1,2,3$.

Since $\lim_{t \to \infty} \omega_e(t)=0$,  by \eqref{peimde1} and \eqref{matrixrhoR}, 
we have
\begin{equation}\label{pe11}
\begin{split}
&\lim_{t \to \infty} (B(t)-\rho_1(t))\\
=&\lim_{t \to \infty} (E_0A(t)-F_1(\omega_e(t))-E_0(NL_1(\omega_e(t))+\zeta(t)))\\
=&\lim_{t \to \infty}(E_0(A(t)-\zeta(t))-F_1(\omega_e(t))-E_0NL_1(\omega_e(t)))\\
=&\lim_{t \to \infty} (E_0(A(t)-\zeta(t)))=0 
\end{split}
\end{equation}
and
\begin{equation}\label{pe12}
\begin{split}
&\lim_{t \to \infty} (E \circ A(t)-\rho_2(t))\\
=&\lim_{t \to \infty} (E \circ A(t)-E \circ (\zeta(t)+NL_1(\omega_e(t))))\\
=&\lim_{t \to \infty} (E \circ(A(t)-\zeta(t))-E \circ NL_1(\omega_e(t)))\\
=&\lim_{t \to \infty} (E \circ(A(t)-\zeta(t)))=0.
\end{split}
\end{equation}
Furthermore, by \eqref{rho3y3}, 
\begin{equation}\label{pe13}
\begin{split}
&\lim_{t \to \infty} E \circ\left(T(\sigma) \varrho(t)-A(t) \mu\right) - E \circ (NL_0(\omega_e(t))-v(t))\\&=\lim_{t \to \infty} E \circ(T(\sigma) \varrho(t)-A(t) \mu+v(t))- E \circ NL_0(\omega_e(t))\\&=0.
\end{split}
\end{equation}

Therefore, by \eqref{equationy12}, \eqref{matrixrhoR},  \eqref{pe11}, \eqref{pe12} and \eqref{pe13}, we have
\begin{equation}\label{pecondition}
\lim_{t \to \infty} (y(t)-\rho (\omega_e(t),\zeta(t),v(t)))=0.
\end{equation}

By \eqref{pecondition0}, we have
\begin{equation}\label{peimconcl1}
\lim _{t \rightarrow \infty} y(t)\tilde{R}(t)=0.
\end{equation}

Since $y(t)$ is PE, there exists a vector $b$ such that the vector $y^T(t)b$ is PE. By \eqref{peimconcl1}, $\lim_{t \to \infty} (y^T(t)b)^T\tilde{R}(t)=0$. From \eqref{adplaw1}, $\lim_{t \to \infty}\dot{\tilde{R}}(t)=\lim_{t \to \infty}\dot{\hat{R}}(t)=0$. Since $y^T(t)b$ is PE, by lemma 4.1 of \cite{liu2009parameter},  we can conclude that $\lim_{t \to \infty}\tilde{R}(t)=0$. The proof is thus completed. 
\end{Proof}


\section{Example}\label{sae}
In this section, the effectiveness of the proposed control approach is demonstrated by its application to the rest-to-rest manoeuvre for a flexible spacecraft system composed of \eqref{system1} and \eqref{system2}, where
the inertia matrix is $J=[J_{11},3,0;3,100,0;0,0,10]$.
which contains an unknown parameter $J_{11}=\mu$, and the damping matrix, the stiffness matrix and the coupling matrix are as follows \cite{ding2014nonsmooth}, where
$C=\rm{diag}(0.1229,0.2195,0.2646,0.1145)$, $K=\rm{diag}(1.2041,\\1.6284,2.7351,5.2409)$ and $\delta=[1.3523,1.1519,2.2167,1.2364;\\$$1.2784,1.0176,1.5891,-1.6537;2.1530,-1.2724,-0.8324,0.2251]$

The disturbance $d_i(t)=C_i\,{\rm sin}(b_{i}t+l_{i})$, $i=1, 2, 3$, where $b_1=1$, $b_2=0.8$, $b_3=\sigma$ is unknown, $C_1=1$, $C_2=2$, $C_3=6$, and $l_1=l_2=l_3=0$.


By \eqref{matrixPhi}, we have $\Phi_{i}=\left[\begin{array}{cc}0 & 1 \\ -b_i^2 &0\end{array}\right], \Psi_{i}=\left[\begin{array}{cc}1 & 0\end{array}\right],~ i=1,2,3.$
Choose $M_{i} = \left[\begin{array}{cc}0 & 1 \\ -3 & -2\end{array}\right], N_{i}=\left[\begin{array}{c}0 \\ 1\end{array}\right],~ i=1,2,3.$
Solve the Sylvester equation yields $T_{i}^{-1}(\sigma)=\left[\begin{array}{cc}3-b_i^2 & 2b_i^2 \\ -2 & 3-b_i^2\end{array}\right],
\Psi_iT_i^{-1}(\sigma)=\left[\begin{array}{ll}
3-b_i^2 & 2
\end{array}\right].$

We can design the internal model of the form \eqref{internal}, where 
$v \in \mathbb{R}^6$, $M=\text{block}\operatorname{diag}\left(M_{1}, M_{2}, M_{3}\right)$, $N=\text{block}\operatorname{diag}\left(N_{1}, N_{2}, N_{3}\right)$, $\Psi=\text{block}\operatorname{diag}\left(\Psi_{1}, \Psi_{2}, \Psi_{3}\right)$, and $T(\sigma)=\text{block}\operatorname{diag}\left(T_{1}(\sigma), T_{2}(\sigma), T_{3}(\sigma)\right)$. Furthermore, we can calculate
\begin{equation}
    L_0(\omega_e)=\left[\begin{array}{ccc}
       0 & 3 & 0\\
       3 & 100 & 0\\
       0 & 0& 10
\end{array}\right] \omega_e,~
    L_1(\omega_e)=\left[\begin{array}{ccc}
       1 & 0 & 0\\
       0 & 0 & 0\\
       0 & 0& 0
\end{array}\right] \omega_e
\end{equation}
and hence $F_0(\omega)$ and $F_1(\omega)$ from \eqref{lpm1}.


Suppose the nominal value of $\sigma=0$, which is the frequency of the disturbance $d_3(t)$. Then, we get
$\Psi T^{-1}(\sigma)=E_0+E^\sigma=E_0+E \Omega(\sigma), \Omega(\sigma)=\sigma^2$
where 
\begin{equation}\label{exsyspar5} 
E_0=\left[\begin{array}{llllll}
2 & 2 & & & & \\
& & 2.36 & 2 & & \\
& & & & 3 & 2
\end{array}\right], E=\left[\begin{array}{lllll}
0 & 0 & & & \\
& & 0 & 0 & \\
& & & -1 & 0
\end{array}\right].
\end{equation}

We can further design
the dynamic compensator \eqref{dyncomp1}, where 
$\zeta \in \mathbb{R}^{6}$ and all parameters are well defined in preceding parameters.


Furthermore, the variables  in \eqref{matrixrhoR} are obtained  as follows:
\begin{equation}
\begin{split}
    \rho\left(\omega_e, \zeta, v\right)=&\left[\rho_1, \rho_2, \rho_3 \right],
    R(\sigma,\mu)=\operatorname{col} (\mu, \sigma^2 \mu, \sigma^2 )    
\end{split}
\end{equation}
where 
$\rho_1=\left[\begin{array}{ccc}
2 \zeta_1+2 \omega_{e 2}&
2.36 \zeta_3+2 \zeta_4 &
3 \zeta_5+2\zeta_6
\end{array}\right]^T+F_1(\omega_e)$,
$\rho_2=\left[\begin{array}{ccc}
0 &
0 &
-\zeta_5
\end{array}\right]^T$, and $\rho_3=\left[\begin{array}{ccc}
0 &
0 &
\eta_5
\end{array}\right]^T$.

By a direct calculation, we can construct an adaptive control law  \eqref{ctrf1}, where $k_1=10$ and $k_2=50$.

\begin{figure}[htp]
\begin{center}
\scalebox{0.55}{\includegraphics{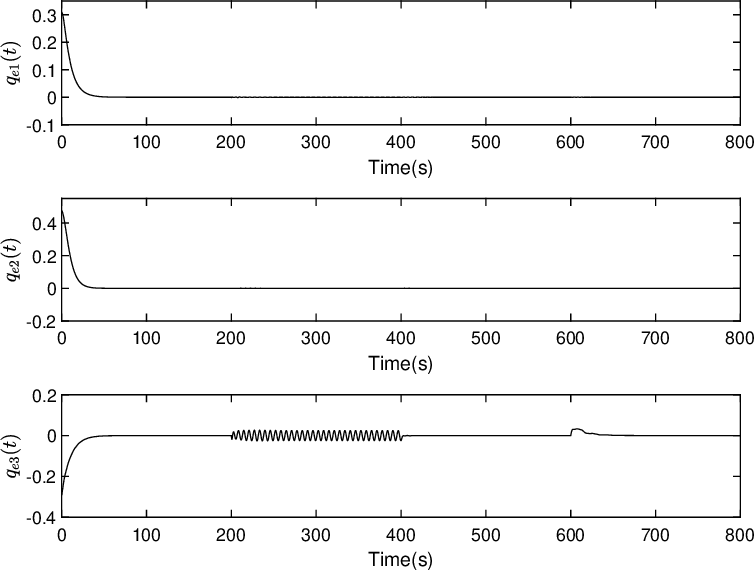}}\caption{ Profile of attitude quaternion errors }\label{fig:qierror}
\end{center}
\end{figure}
\begin{figure}[htp]
\begin{center}
\scalebox{0.55}{\includegraphics{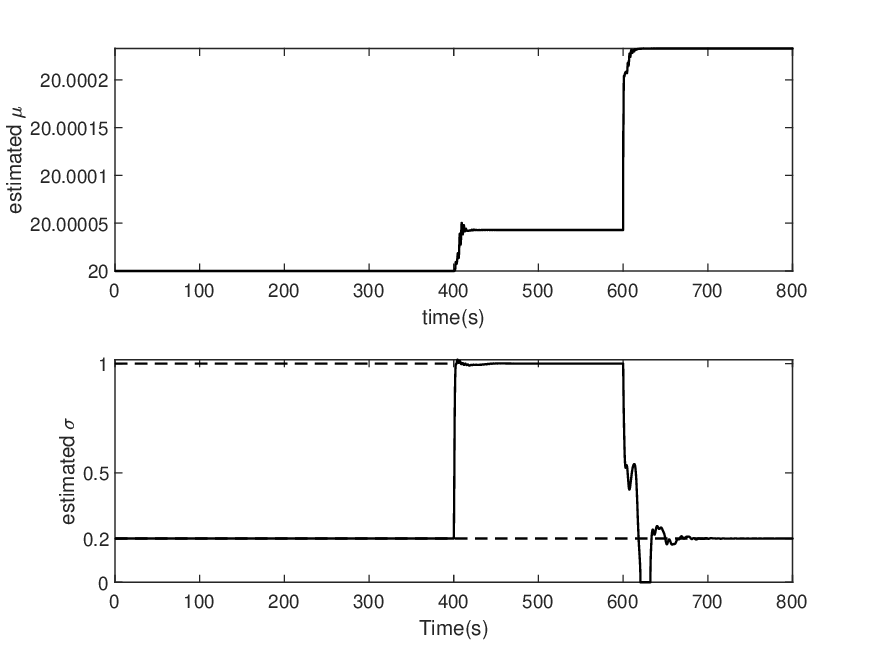}}\caption{Profile of adaptive parameters}\label{fig:estimated}
\end{center}
\end{figure}

During the rest-to-rest manoeuvre, assume that the initial attitude of the flexible spacecraft is
$q(0)=[0.3,-0.2,-0.3,0.8832]^T$, the initial angular velocity is $\omega(0)=[0,0,0]^T$, the desired attitude is $q_d=[-0.24,-0.57,-0.18,0.77]^T$, and the desired angular velocity is $\omega_d=[0,0,0]^T$. The other states' initial values are simply chosen as 0. The simulation results are illustrated in Figs. \ref{fig:qierror}-\ref{fig:estimated}. 
In particular, Fig. \ref{fig:qierror} depicts the tracking performances of $q_i$. During the first 200s, the frequencies of the disturbance and inertia matrix are known,  $\sigma=0.2$, $\mu=20$, and the control law without the update adaptive control law \eqref{adplaw1} is applied.
It is clear that the tracking objective is achieved by the control law. At 200s, the frequency of the disturbance is set to $\sigma=1$ and the inertia matrix parameter is set to $\mu=22$. Consequently, the tracking behaviour is disrupted. At 400s, the adaptive control law is activated with $k=10$. From the simulation result, it can be observed that when the control law is not designed according to the true values, the corresponding control law is still capable of achieving attitude tracking, even when the parameters are reverted to $\sigma=0.2$ and $\mu=20$ at 600s. In particular, for $\hat{R}(\sigma,\mu)$, $\hat{R}_1$ and $\sqrt{\hat{R}_3}$ are the estimated values of $\mu$ and $\sigma$, respectively. From Fig. \ref{fig:estimated}, it can be observed that when the adaptive law is activated, $\sqrt{\hat{R}_3}$ can converge to its true value $\sigma$, but $\hat{R}_1$ cannot converge to  its true value $\mu$. 


Let $A_0=0$, a direct calculation gives $y(t)$ as follows: $y(t)=\rm{diag}(0,0,\frac{\left[\sigma^2-3,2\right] \operatorname{\varrho}_3(t)}{\left(3-\sigma^2\right)^2+4 \sigma^2})$
where $\varrho_3(t)=[d_3(t), \dot{d}_3(t)]^T$.
For $y(t)$ in Theorem \ref{theoremofpe}, define $y(t)=[y_1(t), y_2(t), y_3(t)]$, where $y_1(t)=B(t)$, $y_2(t)=E \circ A(t)$, $y_3(t)=E \circ\left(T(\sigma) \varrho(t)-A(t) \mu\right)$,
we can verify that
\begin{equation}\label{Y3PEequation}
    \lim _{t \rightarrow \infty} y_3(t)\tilde{R}_3(t)=0.
\end{equation}
Let $W(t)=y_3(t)=E \circ[T^{\sigma}\varrho(t)-A(t)\mu]$, we can calculate:
\begin{equation}
    \begin{split}
       W(t)=\left[\begin{array}{ccc}
0 & 0 &
\frac{\left|\sigma^2-3,2\right| \operatorname{\varrho}_3(t)}{\left(3-\sigma^2\right)^2+4 \sigma^2} \\
\end{array}\right]^T.
    \end{split}
\end{equation}
For $\sigma=1$ and $\sigma=0.2$, we can calculate $\varrho_3(t)=\left[\begin{array}{ll}
6sin(t) & 6cos(t)
\end{array}\right]^T$ and $\varrho_3(t)=\left[\begin{array}{ll}
6sin(0.2t) & 3cos(0.2t)
\end{array}\right]^T$, respectively.
Then, for $\sigma=1$, $W(t)=\left[\begin{array}{ccc}
0 & 0 & \frac{3}{2}(cos(t)-sin(t)) 
\end{array}\right]$. Obviously, $W(t)$
is a periodic function with the period $T=2 \pi$. Thus, for any unit vector $c \in \mathbb{R}^2$, $c^T W(t)$ is also a periodic function with period $T=2\pi$, which implies that there exists an $\vartheta_0 \geq 0$ such that
\begin{equation}
    \left(1 / 2\pi \right) \int_t^{t+2\pi}\left|c^{\top} W(s)\right| d s \geqslant \vartheta_0, ~~\forall t \geqslant t_0.
\end{equation}

Furthermore, note that $c \in \mathbb{U}_x=\left\{x \in \mathbb{R}^2,\|x\|=1\right\}$ and $\mathrm{U}_x$ is compact. Therefore, there must exist a minimum value $\vartheta_m>0$ such that for any $c \in \mathbb{U}_x$,
\begin{equation}
\frac{1}{T} \int_t^{t+T}\left|c^T W(s)\right| d s \geq \vartheta_m, \forall t \geq t_0 .    
\end{equation}

It means that $W(t)$ is PE for $\sigma=1$. Similarly, we can show that when $\sigma=0.2$, $W(t)$ is PE. So,  \eqref{Y3PEequation} can be written as follows:
\begin{equation}
\lim _{t \rightarrow \infty}\left[\sigma^2-3,2\right] \varrho_3(t)\left[\hat{R}_3(t)-\sigma^2\right]=0.
\end{equation}
Therefore, we have
$\lim _{t \rightarrow \infty}\hat{R}_3(t)=\sigma^2$,
which is also verified by the simulation result. 



\section{Conclusion}\label{con}
This paper has studied the attitude tracking problem of uncertain flexible spacecraft systems subject to a large class of unknown external disturbances. 
First, a class of nonlinear internal models has been developed, which converts the attitude tracking problem of the uncertain flexible spacecraft system into a regulation problem of an augmented system composed of the flexible spacecraft system and the designed internal model. 
Then,  an auxiliary dynamic system and a dynamic compensator have been developed for auxiliary stabilization and compensation of partial system uncertainty. 
By making a linearly parameterized manipulation and introducing some input transformation, an adaptive control law has been proposed to solve the regulation problem of the augmented system, leading to the solution to the attitude tracking problem under unknown disturbances. In addition, the convergence condition for the estimated parameters has been analyzed based on some PE condition. 
\section*{References}

\bibliographystyle{IEEEtran}

\bibliography{IEEEabrv,mybib_Eul2}

\end{CJK}
\end{document}